\documentclass[11pt]{article}

\usepackage{EAH}
\usepackage{fullpage}
\usepackage{bbm}

\newcommand{\oo}{\mathcal{O}}
\newcommand{\pr}{\mathbb{P}}
\newcommand{\M}{\mathbb{M}}
\newcommand{\I}{\mathbb{I}}
\newcommand{\one}{\mathbbm{1}}
\newcommand{\nrel}{\overline{N\text{-}REL}}
\newcommand{\rhat}{\widehat{R}}

\newcommand{\SET}{\text{SET}}
\newcommand{\BOOL}{\text{BOOL}}
\newcommand{\BIREL}{\text{BI-REL}}
\newcommand{\NREL}{\text{N-REL}}
\newcommand{\kREL}{k\text{-REL}}
\newcommand{\PORDINAL}{\text{PORDINAL}}
\newcommand{\ORDINAL}{\text{ORDINAL}}
\newcommand{\MEAS}{\text{MEAS}}
\newcommand{\PROB}{\text{PROB}}
\newcommand{\INTERVAL}{\text{INTERVAL}}
\newcommand{\SCALAR}{\text{SCALAR}}
\newcommand{\FACE}{\text{FACE}}
\newcommand{\FVECT}{\text{FVECT}}
\newcommand{\RV}{\text{RV}}
\newcommand{\STO}{\text{STO}}

\title{A Category Theoretical Investigation of the Type Hierarchy \\
for Heterogeneous Sensor Integration\thanks{PNNL-25784}}
\author{Emilie Purvine, Cliff Joslyn, Michael Robinson}

\begin{document}
\maketitle

\section{Introduction}
Consider the case of many sensors, each returning very different types of
data (e.g., a camera returning images, a thermometer returning probability
distributions, a newspaper returning articles, a traffic counter returning
numbers). Additionally we have a set of questions, or variables, that we wish
to use these sensors to inform (e.g., temperature, location, crowd size,
topic). Rather than using one sensor to inform each variable we wish to
integrate these sources of data to get more robust and complete information.
The problem, of course, is how to inform a variable, e.g., crowd size, using
a number, a newspaper article, and an image. How do we integrate these very
different types of information? In \cite{RoM2016} Robinson proposes that
sheaf theory is the canonical answer. Moreover, one of the axioms in
\cite{RoM2016} which makes sheaf theory work for data integration is that all
data sources have the structure of a vector space. Therefore, the motivating
question for everything in this report is ``How do we interpret arbitrary
sensor output as a vector space with the intent to integrate?''

The rest of this report is structured as follows.
First, we present the big picture of transforming raw sensor data into vector
space data in Section \ref{sec:bigpicture}. Then in Section
\ref{sec:el_bn_sh} we define category theoretic elements, bundles, and
sheaves. Section \ref{sec:types} contains category theoretic definitions of
many different data types, and in Section \ref{sec:toPMFVECT} we describe the
transformations of each category to the category of pseudo-metrized finite
vector spaces. Finally, in Section \ref{sec:example} we go through an
example of data integration by putting together all of these concepts. 

\section{Big Picture}\label{sec:bigpicture}
In this section we describe a general three step process to transform raw
sensor data into ``cooked'' vector space data. An example will be given at
the end of this report, in Section \ref{sec:example}, after all of the
machinery is built up in the intermediate sections. We begin with the
following assumptions:
\begin{itemize}
  \item Sensor $S$ returns data of the same format with every reading. For
      example, a camera always returns an image or a newspaper always
      returns an article.
  \item Variable $V$ is informed by sensor $S$.
  \item There is an \emph{analytic}, $f_{S,V}$, which takes in a reading
      from $S$ and outputs information that can be interpreted in the
      context of $V$. For example, if $S$ is a newspaper and $V$ is the
      question ``is there violence?'' then $f_{S,V}$ could take an article
      and return the set of words within the article that indicate
      violence. If $S$ informs multiple variables, $V_1, V_2, \ldots$, then
      there is an analytic on $S$ for each variable.
  \item Variable $V$ has a native data type. For example, crowd size is
      numerical, protest topic is categorical, and ``is there violence?''
      is boolean (True/False).
\end{itemize}

\begin{figure}
\begin{center}
\includegraphics[width=0.4\linewidth]{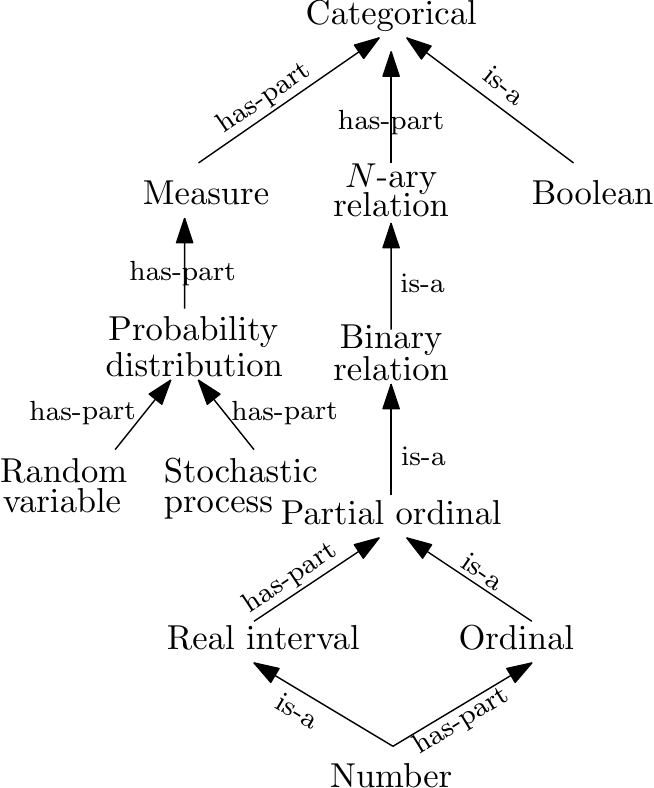}
\caption{Hierarchy of data types returned by analytics.}\label{fig:mathematized}
\end{center}
\end{figure}

Our three step process to interpret output from $S$ as an element of a vector
space begins with the analytic. In Figure \ref{fig:mathematized} we show a
collection of data types that can be returned by analytics. This is not meant
to be exhaustive since other data types certainly exist in the world.
However, it covers the types we need for our purposes. These data types can
be rigorously defined mathematically, and we do so in Section \ref{sec:types}
using the language of \emph{category theory}. At this step it is not required
that we think of the result of the analytic as living in a category, but it
is often easier to do so.

\paragraph{Step 1:}
For output, $s$, from sensor $S$ apply $f_{S,V}$ to \emph{mathematize} the
sensor output. $f_{S,V}(s)$ will be of a type described in Figure
\ref{fig:mathematized}. We can then describe the set of all mathematized data
from sensor $S$ as $f_{S,V}(S)$. This will be a collection of data that all
has the same type. Note that the type returned by the analytic is not
necessarily the same as the native data type of variable $V$.

\vspace{0.25cm} Next our second step moves from the mathematized data into a
category, specifically the category that is native to variable $V$. In Figure
\ref{fig:types} we show the analogous hierarchy to the previous figure, but
now with category names and forgetful, faithful, and inclusion functors.

\paragraph{Step 2:}
Assume variable $V$ has native type $\mathcal{C}$, where $\mathcal{C}$ is a
\emph{category} (to be defined in a later section). Then given the set of
possible results of our analytic, $f_{S,V}(S)$, of some data type from Figure
\ref{fig:mathematized} we identify an object $C \in Ob(\mathcal{C})$. This is
the \emph{cooking} step. We do this in such a way that each $f_{S,V}(s)$ maps
to an \emph{element} of $C$. For a description of elements we refer the
reader to Section \ref{sec:elements}.

\vspace{0.25cm} Our final step in this process is to map each object $C \in
Ob(\mathcal{C})$ to a vector space $W$ so that each element in an object $C$
maps to a single vector $w \in W$. This mapping should be a \emph{functor}
from $\mathcal{C}$ to the category $\FVECT$ of finite dimensional vector
spaces. For the definition of a functor see Section \ref{sec:sheaves}.

\paragraph{Step 3:}
Define a functor from $\mathcal{C}$, the native type for variable $V$, to
\FVECT, the category of finite dimensional vector spaces. If there is
structure to the objects of $\mathcal{C}$ the goal is to reflect that
structure in the image objects in \FVECT.

\begin{figure}
\begin{center}
  \includegraphics[width=0.9\textwidth]{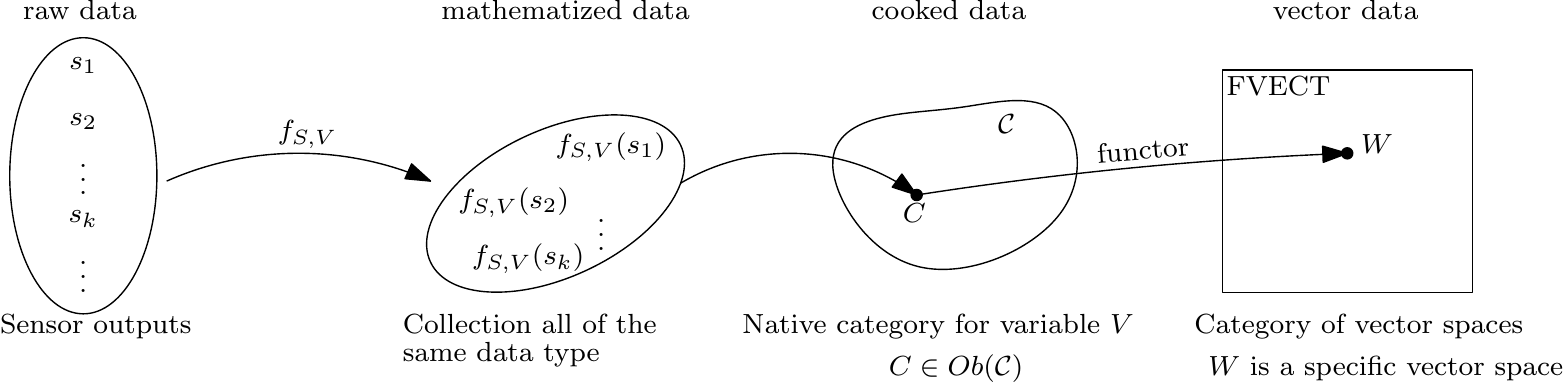}
  \caption{}\label{fig:cooking}
\end{center}
\end{figure}

\vspace{0.25cm} This whole process, which we refer to as
\emph{categorification}, is described pictorially in Figure
\ref{fig:cooking}. We must point out that our use of categorification is
similar to, but distinctly different from other uses of the word. As you read
the remainder of this report, please do so in the context of this pipeline.
In Section \ref{sec:el_bn_sh} it is important to keep in mind that in the
information integration application the stalks will end up being these vector
spaces. While reading Section \ref{sec:types} remember that $\mathcal{C}$
will be one of these categories, and that the sets $f_{S,V}(S)$ can be
interpreted in this context as well. Then, Section \ref{sec:toPMFVECT}
describes the possible functors to \FVECT. Finally, we give an example of
this three step process in Section \ref{sec:example} once we have built up
the terminology to do so.

\section{A note on elements, bundles, sheaves, and
assignments}\label{sec:el_bn_sh} Before going into the type hierarchy in
Section \ref{sec:types} we first formally define categorical elements,
sheaves, and sections. These concepts will be needed as we define our type
hierarchy. We will describe the concepts of a bundle and assignments in
bundles just as defined in Goldblatt \cite{GoR2006}. After that we will do a
similar construction of sheaves and global sections (different from that in
Goldblatt).

It is at this point that we introduce the definition of a category.
\begin{defn}
  A \emph{category}, $\mathcal{C}$, consists of a class of objects,
  $Ob(\mathcal{C})$, and a class of morphisms, $Hom_{\mathcal{C}}(C, D)$, for
  each pair of objects, $C, D \in Ob(\mathcal{C})$. Additionally, for each
  three objects $A, B, C \in Ob(\mathcal{C})$ there is a composition of
  morphisms $\circ :Hom_{\mathcal{C}}(A,B) \times Hom_{\mathcal{C}}(B, C) \rightarrow
  Hom_{\mathcal{C}}(A, C)$, i.e., for each $f:A\rightarrow B$ and
  $g:B\rightarrow C$ there is a unique $g\circ f : A\rightarrow C$, such that
  the following holds:
  \begin{itemize}
  \item Composition of morphisms is associative. If $f:A\rightarrow B$,
      $g:B\rightarrow C$, and $h:C\rightarrow D$ then $h\circ(g\circ f) =
      (h\circ g)\circ f$.
  \item There is an identity morphism for each object C, $1_C:C\rightarrow
      C$, such that $1_C\circ g = g$ and $h\circ 1_C = h$.
  \end{itemize}
\end{defn}

\subsection{Elements}\label{sec:elements}
In category theory, when objects of a category can be arbitrary (i.e., not
small), the concept of an \emph{element} of an object in a category may not
be intuitive. However, we can use morphisms and a terminal object, if one
exists, in the category to define elements of an object. Given a category,
$C$, and two objects, $T, A \in Ob(C)$ we say that morphism $x: T \rightarrow
A$ is a \emph{$T$-valued element of $A$} \cite{BaMWeC2005}. In the case of
small categories, where both $Ob(C)$ and $Hom(C)$ are sets, we do have an
intuitive notion of elements. In order to match our intuition with this
morphism notion of elements we choose $T$ to typically be a terminal object
and call it $\one_{C}$ (these terminal objects will often have size one in
some regard which is why we denote it using $\one$).
\begin{defn}
  A terminal object, $\one$, in a category $C$ is an object such that for any
  other $c \in Ob(C)$ there exists a unique morphism $f \in Hom_C(c, \one)$.
\end{defn}
Then, our elements of $A \in Ob(C)$ will be all morphisms from the chosen
terminal object to $A$. Given this notion of element, in the following
sections we will state what our $\one_C$ object will be in order for us to
choose elements from objects in each category. While it is true that a
terminal object typically defines elements which match our intuition (e.g.,
elements of a set in the category SET), this is not always the case. We will
see in later sections cases in which we use the more general $T$-valued
element of $A$ for some non-terminal object $T$.

\subsection{Bundles and assignments}
A \emph{bundle}, $(\mathcal{A}, p, X)$, is composed of sets $X$ (the
\emph{base space}), and $\mathcal{A}$ (the \emph{stalk space}), and a
function $p : \mathcal{A} \rightarrow X$ mapping elements of the stalk space
to members of the base space. Then for each $x \in X$ the \emph{stalk over
$x$} will be given by $p^{-1}(x) \subset \mathcal{A}$. In other words, the
stalk over $x$ consists of the elements in $\mathcal{A}$ that map to $x$ in
$p$.
\begin{ex}
Let $X = \{1, 2, 3\}$ and $\mathcal{A} = \R \sqcup \Z \sqcup \{1,2,3,4,5\}$
where $\sqcup$ denotes disjoint union (the co-product in \SET), i.e.,
\[\mathcal{A} = \{ \tup{r, 1} : r \in \R \} \cup \{ \tup{z,2} : z \in \Z \}
\cup \{ \tup{i, 3} : 1 \leq i \leq 5, i \in \Z\}.\] Then, define
$p(\tup{x,i}) = i$ which assigns $\R$ to be the stalk over 1, $\Z$ to be the
stalk over 2, and $\{1,2,3,4,5\}$ to be the stalk over 3. Note that this is
fundamentally different than letting $\mathcal{A} = \{ \R, \Z, \{1,2,3,4,5\}
\}$ and $p(\R) = 1$, $p(\Z) = 2$, and $p(\{1,2,3,4,5\}) = 3$.
\end{ex}

Notice that the collection of all bundles over the same base space $X$,
denoted $\text{BN}(X)$, is simply the comma category $\SET\downarrow X$ where
the objects in the category are morphisms in \SET{} that have codomain $X$. A
morphism in $\text{BN}(X)$ from $p_1 : \mathcal{A} \rightarrow X$ to $p_2 :
\mathcal{B} \rightarrow X$ is a morphism $f \in Hom_\SET(\mathcal{A},
\mathcal{B})$ such that the following diagram commutes, i.e., $p_1(a) =
p_2(f(a))$.
\begin{center}
  \includegraphics{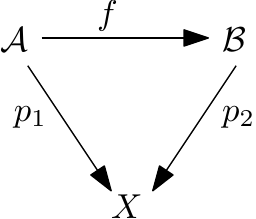}
\end{center}

Informally we have been thinking of assignments as choosing one element from
each stalk in a sheaf (or bundle). We can do this formally using a terminal
object in $\text{BN}(X)$. We claim that the morphism $id_X : X \rightarrow X$
is a terminal object. Consider the following diagram for an arbitrary bundle
$p : \mathcal{A} \rightarrow X$ in $\text{BN}(X)$.
\begin{center}
  \includegraphics{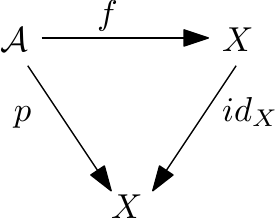}
\end{center}
If this diagram commutes then we have $p(a) = id_x(f(a)) = f(a)$ so $f = p$
is the only choice, making $id_x : X \rightarrow X$ a terminal object in
$\text{BN}(X)$.

Then, how do we understand elements of a particular bundle $p : \mathcal{A}
\rightarrow X$ in $\text{BN}(X)$? Elements are the morphisms $e$ which make
the following diagram commute.
\begin{center}
  \includegraphics{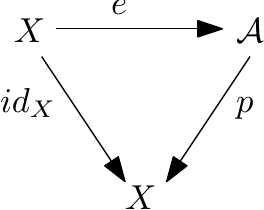}
\end{center}
In other words, $x = p(e(x))$, so $e$ must take each $x$ to an element of its
stalk space in $\mathcal{A}$ as defined by $p$. Essentially then $e$ chooses
one thing from each stalk, which is exactly what we think of as an
assignment.

\subsection{Sheaves}\label{sec:sheaves}

Sheaves are a bit more complicated than bundles. They are more general and
have more structure. But, we want to do something similar to formalize how we
understand assignments. For our purposes, the base space in a sheaf, rather
than being an index set, will be the face category of an abstract simplicial
complex, $X$, which we denote by $\FACE(X)$. In this category the objects are
the faces of $X$ and the morphisms are attachment maps, if $x$ is a subface
of $y$ then $x\rightarrow y$ is a morphism in $\FACE(X)$. Notice that each
morphism is unique which makes $\FACE(X)$ a preorder. Further, since
``subface'' is antisymmetric we know that $\FACE(X)$ is a partial order
category. Then, how do we understand a sheaf? First we must define a
presheaf.
\begin{defn}
  A \emph{presheaf of sets} over an ASC $X$ is a covariant functor, $F : \FACE(X) \rightarrow \SET$, from the face category,
  $\FACE(X)$, to the category \SET.
\end{defn}

\begin{defn}
  A \emph{(covariant) functor}, $F: C\rightarrow D$, from category $C$ to category $D$ is
  a mapping that satisfies the following properties
  \begin{itemize}
    \item For each $X \in Ob(C)$ we have $F(X) \in Ob(D)$
    \item For each morphism $f : X \rightarrow Y$ in $Hom_C(X,Y)$ we map to
        morphism $F(f) : F(X) \rightarrow F(Y)$ in $Hom_D(F(X), F(Y))$ such
        that:
    \begin{itemize}
      \item $F(id_X) = id_{F(X)}$ for every $X \in C$, and
      \item $F(g\circ f) = F(g) \circ F(f)$ for all morphisms $f:
          X\rightarrow Y$ and $g: Y\rightarrow Z$.
    \end{itemize}
  \end{itemize}
\end{defn}

A sheaf is then defined from a presheaf by specifying two axioms which are
called ``locality'' and ``gluing''. We will not go into the specifics here,
but only say that given a presheaf there is a \emph{unique} way of defining a
sheaf. Additionally we point out to the reader that sheaves and presheaves
can be defined in more generality by replacing ASC $X$ with a general
topological space and \SET{} with any concrete category.

The sheaf assigns to each face, $x$, a stalk, $F(x) \in Ob(\SET)$ and to each
attachment map in $\FACE(X)$ a morphism in \SET. This is analogous to a
bundle being a morphism in \SET{} from a stalk space $\mathcal{A}$ to the
base space $X$. Notice that the map goes the other way since we don't want
\emph{every} object in \SET{} to be involved in the sheaf. We have already
observed that morphisms in $\FACE(X)$ are unique. Since it is a category we
have morphism composition so that $x \rightarrow y \rightarrow z$ is equal to
the unique $x\rightarrow z$ which must exist. It is because of this
uniqueness that we guarantee that the resulting morphisms in \SET{} will
commute, i.e., $F(y\rightarrow z)\circ F(x\rightarrow y) = F(w\rightarrow
z)\circ F(x\rightarrow w)$ for all $x\rightarrow y \rightarrow z$ and $x
\rightarrow w \rightarrow z$.

When we introduced bundles we talked about the category of all bundles over a
base space $X$, $\text{BN}(X)$, as being the comma category $\SET\downarrow
X$. So what is the analogous category of all sheaves over the same ASC, $X$?
Let's call it $\text{SH}(X)$. The objects are now functors (instead of
morphisms) from $\FACE(X)$ to \SET, and the morphisms are natural
transformations between functors. This is an example of a \emph{functor
category}.
\begin{align*}
  Ob(\text{SH}(X)) &= \{ F : \FACE(X) \rightarrow \SET ~s.t.~ F \text{ is a functor}\} \\
  Hom_{\text{SH}(X)}(F_1, F_2) &= \{\eta : F_1 \implies F_2 ~s.t.~ \eta \text{ is a natural transformation} \}
\end{align*}
\begin{defn}
  Given two functors, $F_1, F_2 : C \rightarrow D$, from category $C$ to
  category $D$, a \emph{natural transformation}, $\eta : F_1 \implies F_2$,
  has two requirements.
  \begin{enumerate}
    \item To each object $x \in C$ we associate a morphism $\eta_x : F_1(x)
        \rightarrow F_2(x)$ in $D$
    \item For each $f \in Hom_{C}(x,y)$ we must have $\eta_y \circ F_1(f) =
        F_2(f) \circ \eta_x$. In other words, the following diagram must
        commute
    \begin{center}
      \includegraphics{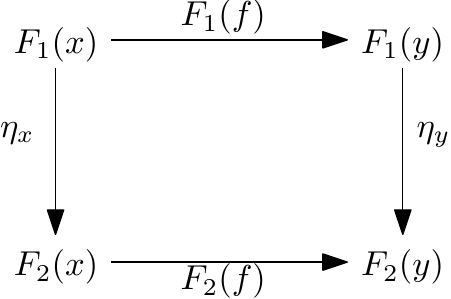}
    \end{center}
  \end{enumerate}
\end{defn}

Recall that we defined assignments in bundles by taking all morphisms from a
terminal object in $\text{BN}(X)$ to our target bundle $p : \mathcal{A}
\rightarrow X$. We can do the same construction as in $\text{SH}(X)$, but we
first need to find the terminal objects in $\text{SH}(X)$. A terminal object
in $\text{SH}(X)$ is a functor, $\mathbbm{1} : \FACE(X) \rightarrow \SET$,
such that there is a unique natural transformation from any other functor,
$F:\FACE(X)\rightarrow \SET$, to $\mathbbm{1}$. If we try to construct a
natural transformation $\eta : F \implies \mathbbm{1}$ and make sure it is
unique we quickly see what it needs to be. For each object $z \in
Ob(\FACE(X))$ we must have a unique $\eta_z : F(z) \rightarrow
\mathbbm{1}(z)$. This means that $\mathbbm{1}(z)$ has to be a terminal object
in \SET. Let $\mathbbm{1}(z) = \{0\}$ for all $z \in \FACE(Z)$. Then for
morphism $f \in Hom_{\FACE(X)}(x,y)$ we define $\mathbbm{1}(f) = id_{\{0\}}$.
One can check that this makes the natural transformation diagram above
commute. So the terminal object in $\text{SH}(X)$ is the functor which sends
every object in $\FACE(X)$ to a terminal object in \SET.

Let's see what happens when we investigate all natural transformations from
our terminal object $\mathbbm{1}$ to another sheaf (functor) $F:\FACE(X)
\rightarrow \SET$. In the previous section this defined assignments of a
bundle by picking out a single element from each stalk. It will do a similar
thing here, but with  more restrictions. Let $\eta : \mathbbm{1} \implies F$.
For each object $z \in \FACE(X)$ we have $\eta_z : \mathbbm{1}(z) \rightarrow
F(z)$. Since $\mathbbm{1}(z)=\{0\}$ a terminal object in \SET, this is the
equivalent of choosing one element from each stalk. But, since we are working
with natural transformations there is one more criteria. For each
$f:x\rightarrow y$ in $Hom_{\FACE(X)}(x,y)$ we must have the following
diagram commute.
\begin{center}
  \includegraphics{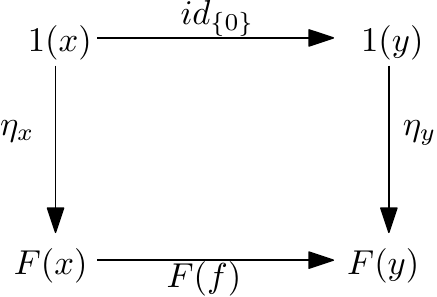}
\end{center}
In words, whatever element that $\eta_x$ chose from stalk $F(x)$ must map,
via $F(f)$, to the element that $\eta_y$ chose from $F(y)$, or $F(f)\circ
\eta_x = \eta_y$. We mentioned that this is stronger than the assignment
criteria, and indeed it is. This defines \emph{global sections}, or
assignments that are globally consistent with respect to the morphisms in the
image $F(\FACE(X)) \subset \SET$. If we want arbitrary assignments we can use
infranatural transformations from $\mathbbm{1}$ which only have property (1)
above.

We built this all up assuming that \SET{} is the target category of the
sheaf. But in fact we could do the same with an arbitrary concrete category
$C$. We can define sheaves of $C$-objects over an ASC $X$, denoted
$\text{SH}_C(X)$, such that the objects are all functors from $\FACE(X)$ to
$C$ and the morphisms are natural transformations. A terminal object in
$\text{SH}_C(X)$ is a functor which sends each $x \in Ob(\FACE(X))$ to a
terminal object in $C$, if one exists, and global sections of a sheaf $F$ are
the natural transformations from a terminal object to $F$.


\section{Type Hierarchy}\label{sec:types}

\begin{figure}[h]
  \centering
  \includegraphics[scale=0.85]{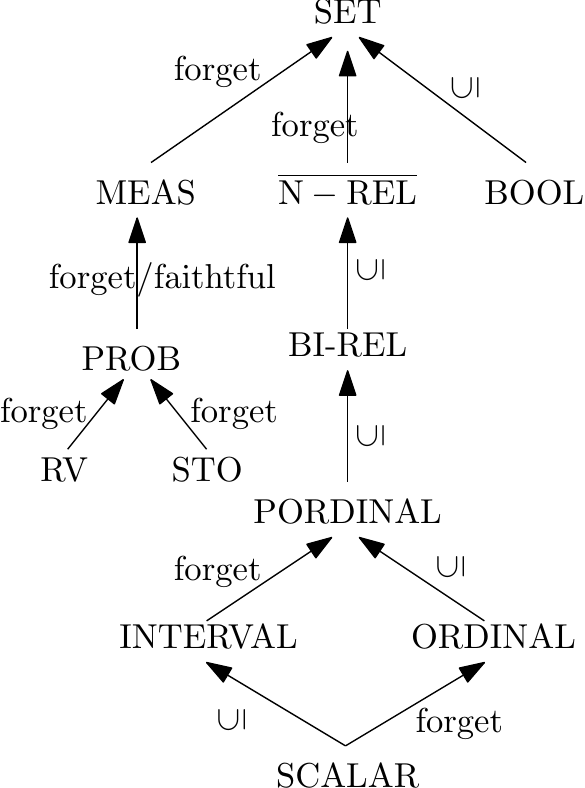}
  \caption{The hierarchy for data types we will consider.}\label{fig:types}
\end{figure}

In this section we define categories for different data types. Data types
that we consider, shown in Figure \ref{fig:types} along with their
relationships, are: categorical, boolean, binary relations, $N$-ary
relations, partial ordinal, ordinal, probability distribution, measurable
spaces, interval-valued, scalar-valued, random variables, and stochastic
processes. For each category we define we will state the \emph{objects} and
\emph{morphisms}, and show that the properties of morphisms are satisfied.
Then, with sheaf theory in mind we state what a stalk would be and finally,
what an assignment would be.

\subsection{Categorical data types, category \SET}
\begin{itemize}
  \item $Ob(\SET)=$ class of all sets
  \item $Hom_{\SET}(S_1, S_2)=$ all set maps from $S_1$ to $S_2$, no
      additional restrictions
      \begin{itemize}
        \item \emph{Composition:} For $f \in Hom_{\SET}(S_1, S_2)$ and $g
            \in Hom_{\SET}(S_2, S_3)$ the composition $g \circ f \in
            Hom_{\SET}(S_1, S_3)$ is the composition of set maps. For $s
            \in S_1$ we have $g\circ f(s) := g(f(s))$.
        \item \emph{Identity maps:} $id_{S} \in Hom_{\SET}(S, S)$ is
            defined to be $id_S(s) = s$ for all $s \in S$.
        \item \emph{Associativity:} For $f: S_1\rightarrow S_2$,
            $g:S_2\rightarrow S_3$, and $h:S_3\rightarrow S_4$ we need to
            show that $h\circ (g\circ f) = (h\circ g)\circ f$. This is
            true since set maps are associative. Both are equal to
            $h(g(f(s)))$ for $s \in S_1$.
      \end{itemize}
  \item \emph{Stalk:} A stalk from \SET{} is a single set $S \in Ob(\SET)$.
  \item \emph{Assignment:} In order to pick elements from an object $S$ to
      make an assignment we use a terminal object $\one_{\SET}=\{0\}$ and
      define elements in $S$ as $Hom_{\SET}(\{0\}, S)$.
\end{itemize}

\subsection{Boolean data types, category \BOOL}
\begin{itemize}
  \item $Ob(\BOOL) = \{\emptyset, \{0\}, \{1\}, \{0,1\} \}$
  \item $Hom_{\BOOL}(B_1, B_2) = Hom_{\SET}(B_1, B_2)$. Composition,
      Identities, and Associativity are inherited from $Hom_{\SET}$.
  \item \emph{Stalk:} A stalk from \BOOL{} is a single set $B \in
      Ob(\BOOL)$.
  \item \emph{Assignment:} We choose a terminal object
      $\one_{\BOOL}=\one_\SET$, so elements are defined as in \SET.
\end{itemize}

\subsection{Binary relation data types, category \BIREL}
Note that this is not the standard category definition for \BIREL.
\begin{itemize}
  \item $Ob(\BIREL) = \{ (S, R) : S \in Ob(\SET), R \subseteq S \times S\}$
      Notice that objects have two parts, a base set and a binary relation.
  \item $Hom_{\BIREL}((S_1, R_1), (S_2, R_2)) =\{ m \in Hom_{SET}(S_1,
      S_2): (x,y) \in R_1 \implies (m(x), m(y)) \in R_2 \}$
      \begin{itemize}
        \item \emph{Composition:} For $m \in Hom_{\BIREL}((S_1, R_1),
            (S_2, R_2))$ and $n \in Hom_{\BIREL}((S_2, R_2), (S_3, R_3))$
            we define $n\circ m$ to be set map composition since
            $m:S_1\rightarrow S_2$ and $n:S_2\rightarrow S_3$. However,
            we must show that $n\circ m$ is indeed a morphism in
            $Hom_{\BIREL}((S_1, R_1), (S_3, R_3))$, i.e., is it relation
            preserving. This is easily seen since both $m$ and $n$ are
            relation preserving. Assume $(x,y) \in R_1$, this implies
            that $(m(x), m(y)) \in R_2$. Then since $n$ is order
            preserving we know that $((n(m(x)), n(m(y)))=(n\circ m(x),
            n\circ m(y)) \in R_3$.
        \item \emph{Identity maps:} Since morphisms here are just
            morphisms in SET we have the identity maps inherited from
            $Hom_{\SET}$. Clearly the identity maps in SET are relation
            preserving.
        \item \emph{Associativity:} This property is also inherited from
            $Hom_{\SET}$.
      \end{itemize}
  \item \emph{Stalk:} A stalk from \BIREL{} is a single binary relation,
      $(S,R) \in Ob(\BIREL)$.
  \item \emph{Assignment:} Our terminal object here is
      $\one_{\BIREL}=(\{0\}, \emptyset)$, the empty relation on a terminal
      object in \SET. The elements of $(S, R)$ are then the morphisms in
      $Hom_{\BIREL}\left((\{0\}, \emptyset), (S, R)\right)$.
\end{itemize}
\subsection{$N$-ary relation data types, category $\overline{\NREL}$}
First we define $\kREL$ for any $k$. In particular we get \BIREL{} when
$k=2$.
\begin{itemize}
  \item $Ob(\kREL) = \{ (S, R) : S \in Ob(\SET), R \subseteq S^k \}$ Again
      objects have two parts, a base set and an $k$-ary relation.
  \item $Hom_{\kREL}((S_1, R_1),(S_2, R_2)) = \{ m \in Hom_{SET}(S_1, S_2)
      : (x_1,\ldots, x_k) \in R_1 \implies (m(x_1), \ldots, m(x_k)) \in R_2
      \}$
      \begin{itemize}
        \item \emph{Composition:} The same construction as above in
            \BIREL{} will give composition for $Hom_{\kREL}$.
        \item \emph{Identity maps:} As in \BIREL, identities are
            inherited from $Hom_{\SET}$.
        \item \emph{Associativity:} This is inherited from $Hom_{\SET}$
            as well.
      \end{itemize}
  \item \emph{Stalk:} A stalk from $\kREL$ is a single $k$-ary relation,
      $(S,R) \in Ob(\kREL)$.
  \item \emph{Assignment:} Our terminal object in \kREL{} is
      $\one_{\kREL}=(\{0\}, \emptyset)$, the empty relation on a terminal
      object in \SET. The elements of $(S, R)$ are then the morphisms in
      $Hom_{\kREL}\left((\{0\}, \emptyset), (S, R)\right)$. Notice that
      this is the same terminal object as in \BIREL. For any $k \in \N$,
      the empty relation is an object in all $\kREL$ categories.
\end{itemize}
Next we define $\overline{\NREL}$ which puts all $k$-ary relations for $2\leq
k\leq N$ into a single category.
\begin{itemize}
  \item $Ob(\overline{\NREL}) = \{ (S, R): S \in Ob(\SET), R \subseteq S^k
      \text{ for some } 2 \leq k \leq N \}$
  \item
  \[Hom_{\overline{\NREL}}((S_1, R_1), (S_2, R_2)) = \left\{
      \begin{array}{ll} \emptyset &
      dim(R_1) \neq dim(R_2) \\
      Hom_{k-\text{REL}}((S_1, R_1), (S_2, R_2)) & dim(R_1)=dim(R_2)=k
      \end{array} \right.\] where $dim(R)$ be the number of elements in
      each relation $r \in R$. In other words, if $R \subseteq S^k$
      then $dim(R) = k$. Composition is inherited from
      $Hom_{k-\text{REL}}$, and therefore so are identity maps and
      associativity.
  \item \emph{Stalk:} A stalk from $\overline{\NREL}$ is a single relation.
  \item \emph{Assignment:} Recall that we observed in \kREL{} that our
      terminal object is the same object no matter what $k$ is. This means
      that $\one_{\overline{\NREL}}=\{\{0\},\emptyset\}$ is a terminal
      object here in $\overline{\NREL}$, and we will pick out the same
      elements from a given object $(S, R)$ regardless of if we are in
      $\kREL$ for a specific $k$ or $\overline{\NREL}$.
\end{itemize}

\subsection{Partial ordinal data types, category \PORDINAL}
\begin{itemize}
  \item $Ob(\PORDINAL) = \{ \po=(P, \mathcal{L}): P \in Ob(\SET),
      \mathcal{L} \subseteq P\times P \text{ is reflexive, transitive, and
      antisymmetric}\}$. Notice that these are binary relations with
      additional properties. Each of these binary relations induces a
      partial order, $\leq$, on $P$ where $p_1 \leq p_2$ iff $(p_1, p_2)
      \in \mathcal{L}$.
  \item $Hom_{\PORDINAL}(\po_1, \po_2) = Hom_{\BIREL}(\po_1, \po_2)$.
      Composition, identity maps, and associativity are inherited from
      \BIREL.
  \item \emph{Stalk:} A stalk from \PORDINAL{} is a single partial order,
      $\po \in Ob(\PORDINAL)$.
  \item \emph{Assignment:} Recall that in \BIREL{} we have $\one_{\BIREL} =
      (\{0\}, \emptyset)$ as our terminal object. Here we cannot choose
      the same object because it does not exist in \PORDINAL. All partial
      orders must be reflexive so we instead have $\one_{\PORDINAL} =
      (\{0\}, (0, 0))$ as a terminal object.
\end{itemize}

\subsection{Ordinal data types, category \ORDINAL}
\begin{itemize}
  \item $Ob(\ORDINAL) = \{ \oo=(O, \mathcal{T}) : O \in Ob(\SET),
      \mathcal{T} \subseteq O\times O \text{ is transitive, antisymmetric,
      and total}\}$. Notice that these are again binary relations with
      additional properties. Each of these induces a total order, $\leq$,
      on $O$ where $o_1 \leq o_2$ iff $(o_1, o_2) \in \mathcal{T}$.
  \item $Hom_{\ORDINAL}(\po_1, \po_2) = Hom_{\BIREL}(\po_1, \po_2)$.
      Composition, identity maps, and associativity are inherited from
      \BIREL.
  \item \emph{Stalk:} A stalk from \ORDINAL{} is a single total order.
  \item \emph{Assignment:} A terminal object in \ORDINAL{} is the same
      as that in \PORDINAL, $\one_{\ORDINAL} = (\{0\}, (0, 0))$.
\end{itemize}

\subsection{Interval data types, category \INTERVAL}
The motivation for creating the \INTERVAL{} category is to let stalks be
subsets of real intervals, $\I$, where $\I = \{ [a,b]: a,b \in \R, a\leq
b\}$. For example, $\I_\Z := \{ [a,b] : a,b \in \Z, a \leq b\}$ is the set of
all integer intervals. We care about both the fact that intervals are
\emph{partially ordered} and they have \emph{algebraic structure} (addition
and multiplication). Therefore, we expand our category to include all
partially ordered semi-rings. The subsets of $\I$ that we care about are
included in this category as objects. We point out here that it is not a
single interval which we claim has additive structure (of course not, $[4,7]$
is an interval but it is not closed under addition since $5+5=10 \not\in
[4,7]$). Instead we consider collections of intervals. We can add, subtract,
and multiply intervals, and there are additive and multiplicative identities.
This leads us to the definition of a semiring.

\begin{defn}
A \textbf{semiring} is a set $S$ with two binary relations, $+$ and $\cdot$,
which are called addition and multiplication respectively such that:
\begin{itemize}
  \item $(R, +)$ is a commutative monoid (operation is associative and
      commutative but does not necessarily have inverses) with identity 0.
  \item $(R, \cdot)$ monoid (operation is associative and does not
      necessarily have inverses) with identity 1.
  \item Multiplication (left and right) distributes over addition.
  \item Multiplication by 0 annihilates $R$.
\end{itemize}
A \textbf{partially ordered semiring} is a semiring, $R$, with a partial
order relation $\leq$ on $R$ satisfying: (a) if $a \leq b$ then $a+c \leq
b+c$, and (b) if $a \leq b$ and $0 \leq c$ then $ac \leq bc$ and $ca \leq
cb$.
\end{defn}
\begin{defn}
A \textbf{partially ordered semiring homomorphism} between $R_1$ and $R_2$ is
a function $f : R_1\rightarrow R_2$ such that $f(a+b) = f(a)+f(b)$, $f(ab) =
f(a)f(b)$, and $f(1_{R_1}) = 1_{R_2}$. Additionally, we require $a \leq b
\Rightarrow f(a) \leq f(b)$.
\end{defn}
Given these definitions we can now describe the category \INTERVAL.
\begin{itemize}
  \item $Ob(\INTERVAL) = \{ \text{partially ordered semirings} \}$.
  \item $Hom_\INTERVAL(R_1, R_2) = \{ \text{partially ordered semiring
      homomorphisms from }R_1\text{ to }R_2 \}$
  \begin{itemize}
    \item \emph{Composition:} This is simply semiring homomorphism
        composition. Given $f \in Hom_\INTERVAL(R_2, R_3)$ and $g \in
        Hom_\INTERVAL(R_1, R_2)$ we define $f\circ g \in
        Hom_\INTERVAL(R_1, R_3)$ to take $r \in R_1$ to $f(g(r)) \in
        R_3$. It is left as an exercise to prove that $f\circ g$ is
        indeed a partially ordered semiring homomorphism.
    \item \emph{Identity maps:} Given an $R \in Ob(\INTERVAL)$ we define
        $id_R \in Hom_\INTERVAL(R, R)$ to be $id_R(r) = r$ for all $r \in
        R$.
    \item \emph{Associativity:} Given $f \in Hom_\INTERVAL(R_2, R_3)$, $g
        \in Hom_\INTERVAL(R_1, R_2)$, and $h \in Hom_\INTERVAL(R_0, R_1)$
        we need to show that $(f\circ g)\circ h = f\circ (g\circ h)$.
        Given an element $r \in R_0$ the left-hand side is
        \[ (f\circ g) \circ h(r) = (f\circ g)(h(r)) = f(g(h(r)))\]
        and the right-hand side is
        \[ f\circ (g\circ h)(r) = f( g\circ h(r)) = f(g(h(r))).\]
        As these are equal we have shown associativity of morphisms.
  \end{itemize}
  \item \emph{Stalk:} A stalk from \INTERVAL{} is any partially ordered
      semiring.
  \item \emph{Assignment:} Typically the $\one_\INTERVAL$ would be the
      terminal object in the category. However, the terminal object in
      \INTERVAL{} is a semiring with one element. When we use that to
      define elements of another semiring $R \in Ob(\INTERVAL)$ we only get
      one element since the semiring with one element is a zero object (it
      is both terminal and initial). So instead we let $\one_\INTERVAL$ be
      a finitely generated semiring with one generator. Then we can choose
      arbitrary elements from other semirings depending on where we map the
      single generator.
\end{itemize}


\subsection{Scalar data types, category \SCALAR}
The motivation for the category \SCALAR{} is similar to that of \INTERVAL{}
but we care about subsets of $\R$ rather than $\I$. In this case $\R$ has a
total order (in contrast to $\I$ which is only partially ordered). So, just
as in \INTERVAL{} where objects  are partially ordered semirings, for
\SCALAR{} we have totally ordered semirings.
\begin{defn}
An \textbf{ordered semiring} is a semiring, $R$, with a total order relation
$\leq$ on $R$ satisfying: (a) if $a \leq b$ then $a+c \leq b+c$, and (b) if
$a \leq b$ and $0 \leq c$ then $ac \leq bc$ and $ca \leq cb$.
\end{defn}
\begin{itemize}
  \item $Ob(\SCALAR) = \{\text{ordered semirings}\}$
  \item $Hom_\SCALAR(R_1, R_2) = \{\text{ordered semiring homomorphisms
      from }R_1\text{ to }R_2\}$ See above for composition, identities, and
      associativity.
  \item \emph{Stalk:} A stalk from the category \SCALAR{} is any ordered
      semi-ring.
  \item \emph{Assignment:} Same as above in \INTERVAL, $\one_\SCALAR$ is a
      finitely generated semiring with a single generator.
\end{itemize}

\subsection{Probability distribution data types, category \PROB}
In this section we will define two types of categories of probability
distributions since there are two types of behaviors we want to capture.
Ultimately we will be modeling output from data sensors as objects in
categories, and both of these behaviors could be expected.
\begin{description}
  \item[Probability distributions:] E.g., a Gaussian with $\mu=0$ and
      $\sigma=1$. For example, this type of data could come from a
      thermometer which may have some error and instead of returning an
      exact temperature, instead returns a probability distribution over
      possible temperatures.
  \item[Sequence of trials of a stochastic process:] For example, this type
      of data could be observations from a camera in a casino pointed at a
      game of blackjack. Each data point is boiled down to the value of
      winnings to a particular player based on the payout random variable
      on the probability space of cards flipped from a deck.
\end{description}
\subsubsection{Elements will be probability distributions} We begin by
describing the first scenario above, where we want elements to be probability
distributions. Because probability distributions are a special kind of
measure we first define the category of all measures, \MEAS. Then we define
\PROB{} as a special case of \MEAS. Let $M$ be a $\sigma$-algebra on a set
$X$. Define $\mathbb{M}(X, M) = \{ \func{m}{M}{\overline{\R}} \text{ where
}m\text{ is a }\sigma\text{-finite measure}\}$. Recall that $\overline{\R} =
\R \cup \{-\infty, +\infty\}$.
\begin{itemize}
  \item $Ob(\MEAS) = \{ \M(X, M) : X \in Ob(SET), M\text{ is a
      }\sigma\text{-algebra on }X\}$
  \item $Hom_{\MEAS}(\M(X,M), \M(Y,N)) = \{ F_\mu: \M(X, M) \rightarrow
      \M(Y, N) ~~s.t.~~ \mu \in \M(X\times Y, M\times N) \}$ where
      \[ F_\mu(m) = n \in \M(Y, N) ~~s.t.~~ n(A) = \int_A \int_X \mu(x,y) m(x) dx dy  = \int_A n(y) dy.\]
      Note that $(X\times Y, M\times N)$ is a measurable space where the
      $\sigma$-algebra is \emph{generated by} $M\times N$, but is not just
      those sets in $M\times N$. We must take all countable unions and
      complements of sets in $M\times N$ in order to define a
      $\sigma$-algebra. We will abuse notation and write $\M(X\times Y,
      M\times N)$ when we mean $\M(X\times Y, \mathcal{M}(M\times N))$
      where $\mathcal{M}(M\times N)$ is the smallest $\sigma$-algebra
      containing $M\times N$.

      We need to prove that $n$ is in fact a measure in $\M(Y,N)$. We must
      check two properties. First, we need to show that $n(\emptyset) =0$.
      From the definition of $n$ we have $n(\emptyset) = \int_\emptyset
      \int_X \mu(x,y) m(x) dx dy$. This is trivially 0 as we are
      integrating over $\emptyset$. Next we need to show that for a
      countable disjoint union of sets $\{ E_j \}_{j=1}^\infty \subseteq N$
      we have $n(\cup_{j=1}^\infty E_j) = \sum_{i=1}^\infty n(E_j)$. This
      is true by additivity of the integral:
      \[ n(\cup_{j=1}^\infty E_j) = \int_{\cup_{j=1}^\infty E_j} \int_X \mu(x,y) m(x) dx dy  = \sum_{i=1}^\infty \int_{E_j} \int_X \mu(x,y) m(x) dx dy = \sum_{i=1}^\infty n(E_j) \]
        \begin{itemize}
        \item \emph{Composition:} Consider $F_\mu \in Hom_{MEAS}(\M(X,M),
            \M(Y,N))$ for $\mu \in \M(X\times Y, M\times N)$ and $F_\nu
            \in Hom_{MEAS}(\M(Y,N), \M(Z,L))$ for $\nu \in \M(Y\times Z,
            N\times L)$.
            \begin{align*}
            \begin{array}{rcccc}
              F_\mu &:& \M(X,M) &\rightarrow &\M(Y,N)\\
                    & &   m     &\mapsto     & n
            \end{array} &\text{ s.t. } n(A) = \int_A \int_X \mu(x,y) m(x) dx dy  \\
            \begin{array}{rcccc}
              F_\nu &:& \M(Y,N) &\rightarrow &\M(Z,L)\\
                    & &   n     &\mapsto     & \ell
            \end{array} &\text{ s.t. } \ell(B) = \int_B \int_Y \nu(y,z) n(y) dy dz
            \end{align*}
            We need to find a $\rho \in \M(X\times Z, M\times L)$ such
            that $F_\nu \circ F_\mu = F_\rho$.
            \begin{align*}
              \ell(B) &= \int_B \int_Y \nu(y,z) n(y) dy dz\\
                      &= \int_B \int_Y \nu(y,z) \int_X \mu(x,y) m(x) dx dy dz\\
                      &= \int_B \int_Y \int_X \nu(y,z)\mu(x,y) m(x) dx dy dz\\
                      &= \int_B \int_X m(x) \int_Y  \nu(y,z)\mu(x,y) dy dx dz \text{ (Fubini)}\\
                      &= \int_B \int_X \rho(x,z) m(x) dx dz
            \end{align*}
            where here $\rho(x,z) = \int_Y  \nu(y,z)\mu(x,y) dy \in
            \M(X\times Z, M\times L)$ as required.
        \item \emph{Identity maps:} Given an $\M(X,M) \in Ob(\MEAS)$ we
            need a $\mu \in \M(X\times X, M\times M)$ so that $F_\mu \in
            Hom_{\MEAS}(\M(X,M), \M(X,M))$ has the property that
            $F_\mu(m) = m$ for all $m \in \M(X,M)$. In other words, we
            need a $\mu$ so that
            \[ m(B) = \int_B \int_X \mu(x,x') m(x) dx dx'\]
            for all $B \in M$. Let
            \[ \mu(x,x') = \left\{\begin{array}{ll}
                                  1 & x=x'\\
                                  0 & else.
                                  \end{array}\right.\]
            Then we can compute the integral
            \begin{align*}
              \int_B \int_X \mu(x,x') m(x) dx dx' &= \int_B \int_B \mu(x,x') m(x) dx dx'\\
                    &= \int_{B\times B} \mu(x,x') m(x) d(x \times x')\\
                    &=  \int_B m(x) dx = m(B)
            \end{align*}
            The first equality is true because $\mu(x,x')$ is definitely
        0 if $x \in B^c$ (since $x' \in B$). The third equality is true
        since $\mu(x,x')=1$ only when $x=x'$ and 0 otherwise.
        \item \emph{Associativity:} Here we need to show that $(F_\nu
            \circ F_\mu) \circ F_\rho = F_\nu \circ (F_\mu \circ
            F_\rho)$. Let
            \begin{align*}
            \begin{array}{rcccc}
              F_\rho &:& \M(W,P) &\rightarrow &\M(X,M)\\
                    & &   p     &\mapsto     & m
            \end{array} &\text{ s.t. } m(C) = \int_C \int_W \rho(w,x) p(w) dw dx \\
            \begin{array}{rcccc}
              F_\mu &:& \M(X,M) &\rightarrow &\M(Y,N)\\
                    & &   m     &\mapsto     & n
            \end{array} &\text{ s.t. } n(A) = \int_A \int_X \mu(x,y) m(x) dx dy  \\
            \begin{array}{rcccc}
              F_\nu &:& \M(Y,N) &\rightarrow &\M(Z,L)\\
                    & &   n     &\mapsto     & \ell
            \end{array} &\text{ s.t. } \ell(B) = \int_B \int_Y \nu(y,z) n(y) dy dz
            \end{align*}
            Let's first look at the LHS. We need to work out $(F_\nu
            \circ F_\mu)$ first, but we have done that already above when
            we defined composition.
            \[
            \begin{array}{rcccccc}
            (F_\nu \circ F_\mu) &:& \M(X,M) &\rightarrow &\M(Y,N) &\rightarrow &\M(Z,L) \\
                                & &  m      &\mapsto     & n      &\mapsto     &\ell
            \end{array}\]
            where $\ell(B) = \int_B \int_X \left[ \int_Y \nu(y,z)
            \mu(x,y) dy\right] m(x) dx dz$. Then we need to compose it
            with $F_\rho$,
            \[
            \begin{array}{rcccccc}
            (F_\nu \circ F_\mu) \circ F_\rho &:& \M(W,P) &\rightarrow &\M(X,M) &\rightarrow &\M(Z,L) \\
                                             & &  p      &\mapsto     & m      &\mapsto     &\ell
            \end{array}\]
            When we work out the details we get
            \begin{align*}
            \ell(B) &= \int_B \int_X \left[\int_Y \nu(y,z) \mu(x,y) dy \right] m(x) ~dx~dz\\
                    &= \int_B \int_X \left[\int_Y \nu(y,z) \mu(x,y) dy \right] \int_W \rho(w,x) p(w) ~dw~dx~dz\\
                    &= \int_B \int_X \int_W \int_Y \nu(y,z) \mu(x,y) \rho(w,x) p(w) ~dy~dw~dx~dz \\
                    &= \int_B \int_W \left[ \int_X \int_Y  \nu(y,z) \mu(x,y) \rho(w,x) dy~dx \right] p(w) dw~dz
            \end{align*}
            Next, we compare to the calculations of the RHS, and confirm
            that we get the same measure in $\M(Z,L)$. Again, we need to
            work out $(F_\mu \circ F_\rho)$ first.
            \[
            \begin{array}{rcccccc}
            (F_\mu \circ F_\rho) &:& \M(W,P) &\rightarrow &\M(X,M) &\rightarrow &\M(Y,N) \\
                                 & &  p      &\mapsto     & m      &\mapsto     &n
            \end{array}\]
            where $n(A) = \int_A \int_W \left[ \int_X \mu(x,y) \rho(w,x)
            dx \right] p(w) dw dy$ as defined by morphism composition.
            Next, we compose with $F_\nu$,
            \[
            \begin{array}{rcccccc}
            F_\nu \circ (F_\mu \circ F_\rho) &:& \M(W,P) &\rightarrow &\M(Y,N) &\rightarrow &\M(Z,L) \\
                                             & &  p      &\mapsto     & n      &\mapsto     &\ell
            \end{array}\]
            When we work the details out here we get
            \begin{align*}
            \ell(B) &= \int_B \int_Y \nu(y,z) n(y) dy dz\\
                    &= \int_B \int_Y \nu(y,z) \left[ \int_W  \int_X \mu(x,y) \rho(w,x) dx ~ p(w) dw \right] dy~dz\\
                    &= \int_B \int_Y \int_W \int_X \nu(y,z) \mu(x,y) \rho(w,x) p(w) dx~dw~dy~dz\\
                    &= \int_B \int_W \left[ \int_Y \int_X \nu(y,z) \mu(x,y) \rho(w,x) dx~dy\right] p(w) dw~dz
            \end{align*}
            This is exactly the same as we got on the LHS (the inner
            integrals are equal by Fubini), so composition is
            associative.

      \end{itemize}
  \item \emph{Stalks:} Stalks in \MEAS{} are objects $\M(X,M)$, i.e., the
      set of measures on $X$ with $\sigma$-algebra $M$.
  \item \emph{Assignment:} In this case we will not be using a terminal
      object, because one does not exist. Instead, we define $\one_{\MEAS}
      = \M(\{0\}, \{ \{0\}, \emptyset \}) = [0,\infty]$. Then given an
      $\M(X, M) \in Ob(\MEAS)$ what does $Hom_{\MEAS}(\one_{\MEAS},
      \M(X,M))$ look like? It is the set of all maps $F_\mu$ for $\mu \in
      \M(\{0\}\times X, \{\{0\}, \emptyset \}\times M)$.
      \begin{align*}
      F_\mu(m)(A) &= \int_A \int_{\{0\}} \mu(y,x) m(y) dy dx\\
                  &= \int_A k \mu(0, x) dx = \int_A \tilde{\mu}(x) dx = \tilde{\mu}(A)
      \end{align*}
      where $k$ is the value of $\int_{\{0\}} m(x) dx = m(\{0\})$, and
      $\tilde{\mu}$ is a measure on $X$, $\tilde{\mu}(x) = \mu(0,x)$.
      Therefore, an element chosen by $\one_{\MEAS}$ from $\M(X, M)$ is
      precisely a single measure over that object's measurable space.
\end{itemize}
Next, define $\pr(X, M) = \{ \func{pr}{M}{\overline{\R}} \text{ where
}pr\text{ is a probability measure}\}$.
\begin{itemize}
  \item $Ob(\PROB) = \{ \pr(X, M) : X \in Ob(SET), M\text{ is a
      }\sigma\text{-algebra on }X\}$
  \item $Hom_{\PROB}(\pr(X,M), \pr(Y,N)) = \{ F_\mu: \pr(X, M) \rightarrow
      \pr(Y, N) ~~s.t.~~ \mu \in \pr(X\times Y, M\times N)\}$ and $\mu$ is
      a conditional probability distribution, i.e., $\int \mu(x,y) dy = 1$
      for all $x$. As before we define
      \[ F_\mu(m) = n \in \pr(Y, N) ~~s.t.~~ n(A) = \int_A \int_X \mu(x,y) m(x) dx dy  = \int_A n(y) dy\]
      \begin{itemize}
        \item \emph{Composition:} Given $F_\mu \in Hom_{\PROB}(\pr(X,M),
            \pr(Y,N))$ and $F_\nu \in Hom_\PROB(\pr(Y,N), \pr(Z,L))$ the
            definition of $F_\nu\circ F_\mu$ is exactly the same as the
            definition in \MEAS. 
        \item \emph{Identity maps:} Given that composition is
            well-defined above, we get identity maps for free from the
            definition in MEAS.
        \item \emph{Associativity:} This follows from the proof of
            associativity in \MEAS.
      \end{itemize}
  \item \emph{Stalks:} Stalks in \PROB{} are objects $\pr(X, M)$, or
      collections of probability measures on $X$ with $\sigma$-algebra $M$.
  \item \emph{Assignments:} Our one element analog for \PROB{} has the same
      $X$ and $M$ as that for \MEAS, $\one_{\PROB} = \pr(\{0\}, \{\{0\},
      \emptyset\})$. However, in this case there is only one probability
      measure on $\{0\}$ which assigns $pr(\{0\}) = 1, pr(\emptyset) = 0$.
\end{itemize}

\subsubsection{Elements will be readings from probability distributions}
We can now discuss the case where elements are trials from a stochastic
process.
\begin{defn}
Let $(\Omega, \mathcal{F}, P)$ be a probability space ($\Omega=$ outcomes,
$\mathcal{F}=$ events is a $\sigma$-algebra over $\Omega$, and
$P:\mathcal{F}\rightarrow[0,1]$ is a probability measure), and $(S, \Sigma)$
be a measurable space ($S$ is a set and $\Sigma$ is a $\sigma$-algebra over
$S$). Then a \emph{random variable}, $X:\Omega\rightarrow S$, is an
$(\mathcal{F}, \Sigma)$-measurable function.

An $S$-valued \emph{stochastic process} is a collection of $S$-valued random
variables on $\Omega$, indexed by a totally ordered set, $T$ (think
``time''). I.e., a stochastic process is a collection $\{X_t: t \in T\}$
where each $X_t$ is an $S$-valued random variable on $\Omega$.
\end{defn}

We ultimately need a category, \STO, of stochastic processes, but we will
begin by defining a category, \RV, of random variables. \STO{} will then be a
generalization of that.
\begin{itemize}
  \item $Ob(\RV) = \left\{ \left[(\Omega, \mathcal{F}, P), (S, \Sigma),
      X:\Omega \rightarrow S\right] \right\}$, in other words, each object
      is a single random variable over a given probability space with a
      given state space.
  \item $Hom_\RV(Y, Z) = $ each morphism will be a pair of maps, $\phi_1,
      \phi_2$, with the following properties\
  \begin{itemize}
    \item $\phi_1 : \Omega_Y \rightarrow \Omega_Z$ is an $(\mathcal{F}_Y,
        \mathcal{F}_Z)$-measurable function
    \item $\phi_2 : S_Y \rightarrow S_Z$ is a $(\Sigma_Y,
        \Sigma_Z)$-measurable function
    \item The following diagram commutes
    \begin{center}
    \includegraphics{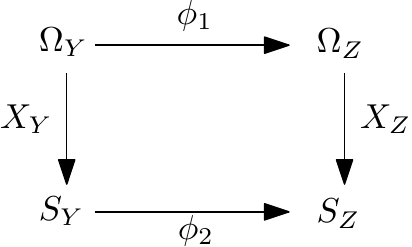}
    \end{center}
    \item \emph{Composition:} Given $(\phi_1, \phi_2) \in Hom_\RV(Y,Z)$,
        and $(\psi_1, \psi_2) \in Hom_\RV(Z,W)$ we define
        $(\psi_1,\psi_2)\circ (\phi_1,\phi_2) \in Hom_\RV(Y,W)$ to be
        $(\psi_1 \circ \phi_1, \psi_2 \circ \phi_2)$. Notice that
        $(\psi_1\circ \phi_1)$ is an $(\mathcal{F}_Y,
        \mathcal{F}_W)$-measurable function. For $F \in \mathcal{F}_W$ we
        have $(\psi_1\circ\phi_1)^{-1}(F) = \phi_1^{-1}(\psi_1^{-1}(F))$,
        and $\psi_1$ is measurable so $\psi_1^{-1}(F)=G$ is measurable in
        $\mathcal{F}_Z$. This implies that $\phi_1^{-1}(G)$ is measurable
        in $\mathcal{F}_Y$. The same argument shows that $(\psi_2\circ
        \phi_2)$ is a $(\Sigma_Y, \Sigma_W)$-measurable function. So
        $(\psi_1 \circ \phi_1, \psi_2 \circ \phi_2)$ is indeed in
        $Hom_\RV(Y,W)$.
    \item \emph{Identity maps:} Given an object $\left[(\Omega,
        \mathcal{F}, P), (S, \Sigma), X:\Omega \rightarrow S\right] \in
        Ob(\RV)$ the identity morphism is $(id_{\Omega}, id_{S})$, the
        identity map on $\Omega$ and the identity map on $S$.
    \item \emph{Associativity:} Function composition is associative so
        composition of these morphisms is also associative.
  \end{itemize}
  \item \emph{Stalk:} A stalk for random variable valued data is a single
      random variable, $\left[(\Omega, \mathcal{F}, P), (S, \Sigma),
        X:\Omega \rightarrow S\right]$.
  \item \emph{Assignment:} First let us ask what a terminal object in \RV{}
      is, and what kind of elements that gives us. A terminal object in
      \RV{} is a random variable $\one_\RV = \left[(\Omega_\one,
      \mathcal{F_\one}, P_\one), (S_\one, \Sigma_\one), X_\one:\Omega_\one
      \rightarrow S_\one \right]$ such that for any other random variable
      $R = \left[(\Omega, \mathcal{F}, P), (S, \Sigma), X:\Omega
      \rightarrow S\right]$ there is a unique morphism
      $R\rightarrow\one_\RV$. So there must be a unique measurable function
      $\Omega \rightarrow \Omega_\one$ and anther unique measurable
      function $S \rightarrow S_\one$. This implies that both $\Omega_\one$
      and $S_\one$ are singleton sets and the $\sigma$-algebras are the
      trivial $\sigma$-algebra containing $\emptyset$ and the set itself.
      This makes a terminal object
      \[ \one_\RV = [(\one_\SET, \{\emptyset,\one_\SET\}, \{P_\one(\emptyset)=0, P(\one_\SET)=1\}), (\one_\SET,\{\emptyset,\one_\SET\}), X_\one : \one_\SET \rightarrow\one_\SET].\]

      Given that this is a terminal object, what does that make our
      elements in \RV? An element of an object $R \in Ob(\RV{})$ would be a
      map from $\one_\RV$ to $R$.
        \begin{align*}
          \phi_1 : \one_\SET \rightarrow \Omega\\
          \phi_2 : \one_\SET \rightarrow S
        \end{align*}
      such that $\phi_1$ and $\phi_2$ are measurable (this is trivial since
      $\phi_i^{-1}(A)$, for $A \subseteq \Omega, S$ is either $\one_\SET$
      or $\emptyset$ and both are measurable), and the corresponding
      diagram commutes. The diagram commuting boils down to the following
      equation
      \[ \phi_2 \circ X_\one(\omega) = X\circ \phi_1 (\omega) \]
      where $\omega$ is the single element of $\Omega_\one$. In other
      words, an element of $R$ is a choice of $x \in \Omega$ and $s \in S$
      such that $X(x) = s$.
\end{itemize}
Given this definition for \RV{} we are now ready to define the category of
stochastic processes for a specific time set, $T$.
\begin{itemize}
  \item $Ob(\STO_T) = \left\{ \left[(\Omega, \mathcal{F}, P), (S, \Sigma),
      \{X_t:\Omega \rightarrow S\}_{t \in T} \right] \right\}$, in other
      words, each object is a stochastic process over a given probability
      space with a given state space.
  \item $Hom_{\STO_T}(Y, Z) = $ each morphism will be two families of maps,
      $\{\phi_{1,t}\}_{t \in T}, \{\phi_{2,t}\}_{t \in T}$, with the
      following properties
    \begin{itemize}
      \item $\phi_{1,t} : \Omega_Y \rightarrow \Omega_Z$ is an
          $(\mathcal{F}_Y, \mathcal{F}_Z)$-measurable function
      \item $\phi_{2,t} : S_Y \rightarrow S_Z$ is a $(\Sigma_Y,
          \Sigma_Z)$-measurable function
      \item The following diagram commutes for all $t \in T$
      \begin{center}
      \includegraphics{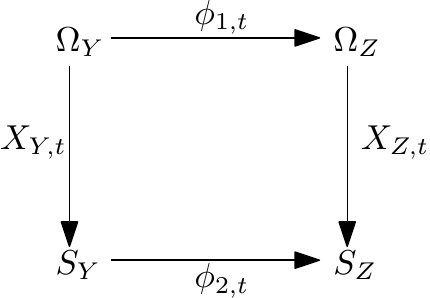}
      \end{center}
      \item \emph{Composition:} The composition of $(\phi_{1,t},
          \phi_{2,t})\in Hom_{\STO_T}(Y,Z)$ and $(\psi_{1,t}, \psi_{2,t})
          \in Hom_{\STO_T}(Z,W)$ is defined just as in \RV{} to be
          \[ (\psi_{1, t} \circ \phi_{1, t}, \psi_{2,t} \circ \phi_{1,t}) \in Hom_{\STO_T}(Y,W). \]
          The compositions are measureable for the same reasons as in the
          category \RV.
      \item \emph{Identity:} Identities are also defined to be just the
          identity functions on $\Omega_Y$ and $S_Y$.
      \item \emph{Associativity:} Again, function composition is
          associative so these morphisms are associative.
    \end{itemize}
    \item \emph{Stalks:} A stalk for a stochastic process indexed by $T$ is
        a single object in $\STO_T$
        \[ \left[(\Omega, \mathcal{F}, P), (S, \Sigma),
      \{X_t:\Omega \rightarrow S\}_{t \in T} \right]. \]
    \item \emph{Assignment:} Just as in the case of \RV{} we need to
        discover the structure of a terminal object in $\STO_T$ in order to
        define elements, and thus assignments. A terminal object in
        $\STO_T$ is a stochastic process
        \[ \one_{\STO_T} = \left[(\Omega_\one,\mathcal{F_\one}, P_\one), (S_\one, \Sigma_\one), \{X_{\one,t}:\Omega_\one \rightarrow S_\one\}_{t \in T} \right]\]
        such that for any other stochastic process $R \in Ob(\STO_T)$ there
        is a unique morphism from $R$ to $\one_{\STO_T}$. Again this
        implies that $\Omega_\one$ and $S_\one$ are terminal objects in
        \SET, and $\mathcal{F}_\one$ and $\Sigma_\one$ are trivial
        $\sigma$-algebras. Then all of the $X_{\one,t}$ are forced to be
        identical. An element in a stochastic process object, $R$, is found
        using maps from $\one_{\STO_T}$ to $R$. This is now two families of
        maps
        \begin{align*}
          \phi_{1,t} : \one_\SET \rightarrow \Omega\\
          \phi_{2,t} : \one_\SET \rightarrow S
        \end{align*}
      such that $\phi_{1,t}$ and $\phi_{2,t}$ are measurable for all $t \in
      T$ (as before this is trivial), and the corresponding diagrams
      commute. The diagrams commuting boil down to the following family of
      equations
      \[ \phi_{2,t} \circ X_{\one,t}(\omega) = X_t\circ \phi_{1,t} (\omega) \]
      where $\omega$ is the single element of $\Omega_\one$. In other
      words, an element of $R$ is a choice of $\{x_t\} \subset \Omega$ and
      $\{s_t\} \subset S$ such that $X_t(x_t) = s_t$ for all $t \in T$.
\end{itemize}

\subsection{Maps in the type hierarchy}
Now that we have defined all of the categories in the type hierarchy we can
fill in the maps between them. For each arrow in the type hierarchy we will
define a functor which describes the transformation.
In the rest of this section as we describe the functors between the
categories it is left as an exercise to show that they respect $F(id_X) =
id_{F(X)}$ and composition as required.

\subsubsection{Inclusion functors}
The following pairs of categories admit an inclusion functor from the first
to the second. In these cases the functors are trivial to define and they
respect $F(id_X) = id_{F(X)}$ and $F(g\circ f) = F(g) \circ F(f)$ as required
in order to be a functor.
\begin{itemize}
  \item $\BOOL \rightarrow \SET$: All objects in \BOOL{} are also objects
      in \SET, and morphisms are the same as in \SET.
  \item $\PORDINAL \rightarrow \BIREL$: Every partial ordered is a binary
      relation, and every order preserving map is a relation preserving map
  \item $\ORDINAL \rightarrow \PORDINAL$: Every total order is a partial
      order
  \item $\kREL \rightarrow \overline{\NREL}$ (in particular, $\BIREL
      \rightarrow \overline{\NREL}$): Clear from  how $\overline{\NREL}$ is
      defined.
  \item $\SCALAR \rightarrow \INTERVAL$: Every ordered semiring is a
      partially ordered semiring.
\end{itemize}

\subsubsection{Non-inclusion functors}
The rest of the functors we will describe are non-trival. Many are still
straightforward, but they are not inclusion maps like those above.

\begin{description}
\item[$F_{SO}: \SCALAR \rightarrow \ORDINAL$] Let $S \in Ob(\SCALAR)$ be an
    ordered semiring with total order $\leq$, and $f \in Hom_\SCALAR(S_1,
    S_2)$ an ordered semiring homomorphism in \SCALAR.
    \begin{itemize}
      \item $F_{SO}(S) = (set(S), \mathcal{T})$. We define $\mathcal{T}$
          to be the total order binary relation induced by $\leq$ where
          $(s,r) \in \mathcal{T}$ iff $s\leq r$ in the ordered semiring.
          Additionally $set(S)$ is the set of elements from the semiring
          $S$, i.e., we forget the semiring structure of $S$.
      \item $F_{SO}(f) = f$. We know that any ordered semiring
          homomorphism is order preserving.
    \end{itemize}
\item[$F_{IP}: \INTERVAL \rightarrow \PORDINAL$] Let $I \in Ob(\INTERVAL)$
    a partially ordered semiring with partial order $\leq$, and $f \in
    Hom_\INTERVAL(I_1, I_2)$ a partially ordered semiring homomorphism.
    \begin{itemize}
      \item $F_{IP}(I) = (set(I), \mathcal{L})$. We define $\mathcal{L}$
          to be the total order binary relation induced by $\leq$ where
          $(i,j) \in \mathcal{L}$ iff $i\leq j$ in the partially ordered
          semiring. Additionally $set(I)$ is the set of elements from the
          semiring $I$, i.e., we forget the semiring structure of $I$.
      \item $F_{IP}(f) = f$. We know that any partially ordered semiring
          homomorphism is order preserving.
    \end{itemize}
\item[$F_{NS} : \nrel \rightarrow \SET$] Let $(S, R) \in Ob(\nrel)$ be an
    object in $\nrel$, and $m \in Hom_{\nrel}((S_1,R_1), (S_2, R_2))$.
    \begin{itemize}
      \item $F_{NS}( (S, R) ) = S$. Recall that $(S,R) \in Ob(\nrel)$ is
          a set $S$ with a $k$-ary relation, $R \subseteq S^k$ for some
          $2\leq k\leq N$. Therefore, we can simply map $(S,R)$ to its
          underlying set $S$ forgetting about the relation structure.
      \item $F_{NS}(m) = m$. Again, recall that any map $m \in
          Hom(\nrel)$ is simply a map $m \in Hom_{\SET}(S_1, S_2)$ with
          extra restrictions (order preserving with respect to $R_1$ and
          $R_2$). So, since $m$ exists in both categories this map is
          allowed.
    \end{itemize}
    This is a forgetful functor. Since the morphism map is inclusion it
automatically respects $F_{NS}(id_{S,R}) = id_{F_{NS}((S,R))}$ and
composition of morphisms.
\item[$F_{PM}: \PROB \rightarrow \MEAS$] Let $\pr(X, M) \in Ob(\PROB)$ be a
    set of probability measures over $X$ with $\sigma$-algebra $M$, and
    $F_\mu \in Hom_\PROB(\pr(X, M), \pr(Y, N )$ a morphism in \PROB.
    \begin{itemize}
      \item $F_{PM}(\pr(X,M)) = \M(X,M)$ this simply expands the
          probability measure space to a generic measure space. Note that
          $\pr(X, M) \subseteq F_{PM}(\pr(X,M))$.
      \item $F_{PM}(F_\mu) = F_\mu$
    \end{itemize}
\item[$F_{MS}: \MEAS \rightarrow \SET$] Let $\M(X, M) \in Ob(\MEAS)$ be the
    set of measures over $X$ with $\sigma$-algebra $M$, and $F_\mu \in
    Hom_\MEAS(\M(X, M), \M(Y, N ))$ a morphism in \MEAS.
    \begin{itemize}
      \item $F_{MS}(\M(X,M)) = \M(X,M)$
      \item $F_{MS}(F_\mu) =$ the map induced by $F_\mu$.
    \end{itemize}
\item[$F_{RP} : \RV \rightarrow \PROB$] Let $\left[(\Omega, \mathcal{F},
    P), (S, \Sigma), X:\Omega \rightarrow S\right] \in Ob(\RV)$ be an
    object in the category \RV, and $(\phi_1,\phi_2) \in Hom_\RV(Y,Z)$ be a
    morphism.
    \begin{itemize}
      \item $F_{RP}(\left[(\Omega, \mathcal{F}, P), (S, \Sigma), X:\Omega
          \rightarrow S\right]) = \pr(S,\Sigma)$. Note that
          $F_{RP}(\one_\RV) = \one_\PROB$ as we would like.
      \item $F_{RP}((\phi_1,\phi_2)) = \delta_{y,\phi_2(y)} \in \pr(S_Y
          \times S_Z, \Sigma_Y\times \Sigma_Z)$. This is the Dirac delta
          function in which $\delta_{y,\phi_2(y)}(s_Y, s_Z) =
          \delta(s_Z-\phi_2(s_Y)) = 0$ unless $s_Z-\phi_2(s_Y)=0$. The
          value on $(s_Z, \phi_2(s_Y))$ is such that
          \[ \int_{S_Y\times S_Z} \delta_{y,\phi_2(y)} (s_Y,s_Z) ds_Y ds_Z = 1. \]
    \end{itemize}
\item[$F_{SP} : \STO \rightarrow \PROB$] Let $\left[(\Omega, \mathcal{F},
    P), (S, \Sigma), \{X_t:\Omega \rightarrow S\}_{t \in T}\right] \in
    Ob(\RV)$ be an object in the category \STO, and $(\{\phi_{1,t}\}_{t \in
    T},\{\phi_{2,t}\}_{t\in T}) \in Hom_\STO(Y,Z)$ be a morphism.
    \begin{itemize}
      \item $F_{SP}(\left[(\Omega, \mathcal{F}, P), (S, \Sigma),
          \{X_t:\Omega \rightarrow S\}_{t \in T}\right]) =
          \pr(S^T,\Sigma^T)$. This is the set of all probability
          distributions over $S\times S\times\cdots$ where there are $T$
          copies of $S$, with $\sigma$-algebra $\Sigma^T$.
      \item $F_{SP}(\{\phi_{1,t}\}_{t \in T},\{\phi_{2,t}\}_{t\in T}) =
          \delta_{\{y_t\}, \{\phi_{2,t}(y_t)\}} \in \pr(S_Y^T\times
          S_Z^T, \Sigma_Y^T\times\Sigma_Z^T)$. This again is the Dirac
          delta function. This time we define
          \[ \delta_{\{y_t\}, \{\phi_{2,t}(y_t)\}}\left(\tup{s_{Y,t}}_{t \in T}, \tup{s_{Z,t}}_{t \in T}\right) =
                \delta\left(\tup{s_{Z,t} - \phi_{2,t}(s_{Y,t})}_{t \in
          T}\right)\] which equals 0 unless $s_{Z,t}-\phi_{2,t}(s_{Y,t})
          = 0$ for all $t\in T$. The other values are such that the total
          integral is 1 just as in the previous case.

    \end{itemize}
\end{description}

%

\section{Mapping to FVECT}\label{sec:toPMFVECT}

We have now reached the point of mapping our data from a category defined in
the previous section to a finite dimensional vector space. As we described in
Section \ref{sec:bigpicture} we need a functor from each category to FVECT,
and our goal is to preserve structure of objects wherever possible.

\subsection{Mapping of objects in \SET}
The objects in \SET{} do not have any structure, they are simply collections
of unique elements. Therefore, we have a relatively simple functor, $F : \SET
\rightarrow \FVECT$, defined as follows:
\begin{itemize}
  \item For $S \in Ob(\SET)$ we define $F(S) = \R[S]$, the
      $|S|$-dimensional vector space with basis being the elements of $S$.
      This can also be written as $\R^S$.
  \item For morphism $f : S_1 \rightarrow S_2$ in $Hom_\SET(S_1, S_2)$ we
      define $F(f) : \R[S_1] \rightarrow \R[S_2]$. Given $v \in \R[S_1]$
      its image $F(f)(v)$ has coefficient of basis element $s \in S_2$
      equal to the sum of coefficients in $v$ from basis elements in
      $f^{-1}(s) \subseteq S_1$. One can easily check that this is a linear
      transformation, so $F(f) \in Hom_\FVECT(\R[S_1], \R[S_2])$, and that
      $F$ satisfies the two requirements of being a functor.
\end{itemize}
We point out here that although each element $s$ of object $S$ picks out a
unique element in $\R[S]$, $1\times s$, we cannot do the reverse. Given $s+t
\in \R[S]$ there is no element in set $S$ which maps to it. In other words,
this functor induces a function from $S$ to $\R[S]$ which is one to one but
not onto.
%

\subsection{Mapping of objects in \BOOL}
\BOOL{} is a subcategory of \SET{} so the categorification is exactly the
same for any $B \in Ob(\BOOL)$.

%

\subsection{Mapping of objects in \kREL}
Let $(S, R) \in Ob(\kREL)$ be a $k$-ary relation, so that $S=\{s_1,\ldots,
s_n\}$ and $R \subseteq S^k$. In this case there is a significant amount of
structure in the $k$-ary relation that we wish to translate into \FVECT.
First we will describe how the functor acts on objects and give some examples
and then we will describe the functor on morprhisms. Consider the vector
space $\R[\rhat] = \left\{\sum_{r \in \rhat} a_r\cdot r : a_r \in \R\right\}$
where $\rhat = R \cup \left\{\tup{s_i,\ldots,s_i}\right\}_{i=1}^n$ is the
relation extended by all reflexive relations. This is isomorphic to
$\R^{|\rhat|}$, an $|\rhat|$-dimensional real vector space and carries
information about all of the relations (just in the names of the basis
elements), but does not tell how they fit together. The vector space that we
will assign to $(S,R)$ is a subspace of $\R[\rhat]$. In particular it is the
subspace spanned by $|S|$ vectors, one for each element of $S$. Notice that
$|S| < |\rhat|$ since we have added the reflexive elements to the relation.
Let $\rhat(s_i) = \{r \in \rhat : s_i \in r\}$ be the set of relations which
$s_i$ is involved in. Given an enumeration of relations $\hat{R} = \{r_1,
r_2,\ldots, r_m\}$ we can consider the column vector
\[ \overrightarrow{\rhat(s_i)} = \tup{a_j}_{j=1}^m, \text{ where } a_j = \left\{
\begin{array}{ll}
1 & r_j \in \rhat(s_i)\\
0 & \text{else}
\end{array}\right..
\]
Then, the subspace of $\R[\rhat]$ that we assign to $(S,R)$ is the space
spanned by $\left\{\overrightarrow{\rhat(s_i)}\right\}_{i=1}^n$.

\paragraph{Example:} Consider the set $S = \{a,b,c,d,e\}$ and ternary relation $R =
\{\tup{a,b,c}, \tup{b,c,e}, \tup{c,a,e}, \tup{d,b,e}\}$. Our base space
$\R[\rhat]$ is defined as
\begin{align*}
\R[\rhat] =& \left\{ \alpha_1\cdot\tup{a,b,c}+\alpha_2\cdot\tup{b,c,e}+\alpha_3\cdot\tup{c,a,e}+\alpha_4\cdot\tup{d,b,e} +\right.\\
           &\left.+\alpha_5\cdot\tup{a,a,a} + \alpha_6\cdot\tup{b,b,b} + \alpha_7\cdot\tup{c,c,c} + \alpha_8\cdot\tup{d,d,d} + \alpha_9\cdot\tup{e,e,e} : \alpha_i \in \R\right\}
\end{align*}
Next, for each $s_i \in S$ we define $\rhat(s_i)$:
\begin{align*}
  \rhat(a) &= \left\{ \tup{a,a,a}, \tup{a,b,c}, \tup{c,a,e} \right\}\\
  \rhat(b) &= \left\{ \tup{b,b,b}, \tup{a,b,c}, \tup{b,c,e}, \tup{d,b,e} \right\}\\
  \rhat(c) &= \left\{ \tup{c,c,c}, \tup{a,b,c}, \tup{b,c,e}, \tup{c,a,e} \right\}\\
  \rhat(d) &= \left\{ \tup{d,d,d}, \tup{d,b,e} \right\}\\
  \rhat(e) &= \left\{ \tup{e,e,e}, \tup{b,c,e}, \tup{c,a,e}, \tup{d,b,e} \right\}.
\end{align*}
Given this, and the ordering of elements of $\rhat$ above in the definition
of $\R[\rhat]$ we can see that the $|S|$-dimensional subspace we want is
defined to be the span of the following five column vectors, $\left\{
\overrightarrow{\rhat(s_i)}\right\}_{i=1}^5 $:
\[
\left[
\begin{array}{c}
1 \\
0 \\
1 \\
0 \\
1 \\
0 \\
0 \\
0 \\
0
\end{array}
\right],
\left[
\begin{array}{c}
1 \\
1 \\
0 \\
1 \\
0 \\
1 \\
0 \\
0 \\
0
\end{array}
\right],
\left[
\begin{array}{c}
1 \\
1 \\
1 \\
0 \\
0 \\
0 \\
1 \\
0 \\
0
\end{array}
\right],
\left[
\begin{array}{c}
0 \\
0 \\
0 \\
1 \\
0 \\
0 \\
0 \\
1 \\
0
\end{array}
\right],
\left[
\begin{array}{c}
0 \\
1 \\
1 \\
1 \\
0 \\
0 \\
0 \\
0 \\
1
\end{array}
\right].
\]
Notice that these are all linearly independent because of the extension of
$R$ to $\rhat$.

This construction of the vector space from $(S,R)$ is the functor on objects
from \kREL{} to \FVECT. Now we need to specify how morphisms in \kREL{} are
mapped to \FVECT. Let $f \in Hom_\kREL((S_1,R_1),(S_2,R_2))$. So $f : S_1
\rightarrow S_2$ is a relation-preserving set function, i.e., if
$\tup{s_{i_1},s_{i_2},\ldots, s_{i_k}} \in R_1$ then $\tup{f(s_{i_1}),
f(s_{i_2}),\ldots, f(s_{i_k})} \in R_2$. We must define $F(f) \in
Hom_\FVECT\left(F(S_1, R_1), F(S_2, R_2)\right)$. Because $F(f)$ must be a
linear transformation it is enough to define the function on the basis
elements of $F(S_1, R_1)$. In both $F(S_1, R_1)$ and $F(S_2, R_2)$ there is a
basis element for each of the elements in $S_1$ and $S_2$ respectively.
Therefore, we can define $F(f)$ by mapping $\overrightarrow{\rhat_1(s_i)}$ to
$\overrightarrow{\rhat_2(f(s_i))}$. That is, map the basis element of $F(S_1,
R_1)$ corresponding to $s_i \in S_1$ to the basis element of $F(S_2, R_2)$
corresponding to $f(s_i) \in S_2$.

\subsection{Mapping of objects in \PORDINAL}
\PORDINAL{} is a subcategory of \BIREL{} and so we can use the
functor described for \kREL.

%
%
%

\subsection{Mapping of objects in \ORDINAL}
Since \ORDINAL{} is a subcategory of \PORDINAL{} the functor is inherited.

\subsection{Mapping of objects in \PROB{} and \MEAS}
The objects in \MEAS{} can be thought of as vector spaces with an extension
to negative measures. First, consider an object $\M(X,M) \in Ob(\MEAS)$. This
consists of all $\sigma$-finite measures, $m:M\rightarrow \overline{\R}$. All
measures are positive, i.e., for all $s \in M$ we have $m(s) \geq 0$.
However, if we additionally allow measures to be totally negative (i.e., for
all $s \in M$ we have $m(s) \leq 0$) we can treat this as a vector space.
Given $m_1,m_2 : M \rightarrow \overline{\R}$ we define $m_1+m_2$ on $M$ by
letting $(m_1+m_2)(s) = m_1(s)+m_2(s)$. Let $\overline{\M}(X,M)$ be $\M(X,M)
\cup -\M(X,M)$, i.e., the set of all measures union the set of all negative
measures. This satisfies the axioms of a vector space:
\begin{itemize}
  \item Associativity and commutativity are clear
  \item The identity element is the function $\mathbf{0} : M\rightarrow
      \overline{R}$ where $\mathbf(0)(s) = 0$  for all $s \in M$.
  \item Additive inverses are simply the negative measure for any given
      measure. For $m \in \M(X,M)$ the negative measure $-m$ is defined to
      be $(-m)(s) = -1\cdot m(s)$ for all $s \in M$.
  \item If for a scalar $a \in \R$ we define scalar multiplication as
      $(am)(s) = a\cdot m(s)$ for all $s \in M$ then this satisfies $a(bm)
      = (ab)m$.
  \item This scalar multiplication is clearly distributive, $a(m_1+m_2) =
      am_1+am_2$ and $(a+b)m = am+bm$.
\end{itemize}
Using this extension we can define our functor $F: \MEAS \rightarrow \FVECT$
as follows:
\begin{itemize}
\item Given object $\M(X,M) \in Ob(\MEAS)$ we have $F(\M(X,M)) =
    \overline{\M}(X,M)$
\item Let $f \in Hom_\MEAS(\M(X, M), \M(Y, N))$, then there is a $F(f) \in
    Hom_\FVECT\left(\overline{\M}(X,M), \overline{\M}(Y,N)\right)$ which
    extends $f$ to the negative measures. This induced map is indeed a
    linear transformation of these vector spaces, a fact which is left up
    to the reader to verify.
\end{itemize}

Notice that unlike \MEAS, \PROB{} is not closed under addition or scalar
multiplication. However, it is closed under convex combination. Given a
collection of probability measures $p_i \in \pr(X,M) \in Ob(\PROB)$ we can
form a new probability measure $\sum a_i p_i$ if all $a_1 \geq 0$ and $\sum
a_i = 1$. So, \PROB{} forms a convex subset of \MEAS. We should be able to
map \PROB{} to \FVECT{} in the same way as \MEAS{}. In other words, we
consider elements of \PROB{} to be elements of \MEAS{}.

\subsection{Mapping of objects in \INTERVAL}
Let $R \in Ob(\INTERVAL)$ so that $R$ is a partially ordered semiring. Since
$R$ has a partial order we could use the \PORDINAL{} functor. But, this does
not take into account the semiring structure of the objects. We will be
continuing to study possible functors for \INTERVAL{} which preserve all
structure within the semirings.

\subsection{Mapping of objects in \SCALAR}
\SCALAR{} is a subcategory of \INTERVAL{} so we define the functor in the
same way.

\subsection{Mapping of objects in \RV{} and \STO}
Let $[(\Omega, \mathcal{F}, P), (S,\Sigma), X:\Omega\rightarrow S] \in
Ob(\RV)$ be a random variable object. Before attempting to create a functor
to \FVECT{} we must first ask, what is the structure that we wish to
preserve? We may wish to preserve the information contained in the $X$ random
variable map. 
There are certainly other kinds of structure in this object that one might
wish to preserve. But in the case of the $X$ random variable map, we can use
the functor we defined for $\SET\rightarrow\FVECT$. Consider $\{(\omega,
X(\omega))\}_{\omega \in \Omega}$ as an object in \SET{} and let $F([(\Omega,
\mathcal{F}, P), (S,\Sigma), X:\Omega\rightarrow S])$ be the image of
$\{(\omega, X(\omega))\}_{\omega \in \Omega}$ under the
$\SET\rightarrow\FVECT$ functor. A similar functor can be constructed for
\STO. We will continue to investigate other possible functors from \RV{}
which preserve other types of structure within the objects.

\section{An example}\label{sec:example}
\begin{figure}[h]
\begin{center}
  \includegraphics{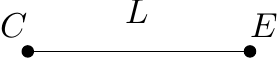}
  \caption{The abstract simplicial complex for the variable $L$ with informed
  by two sensors, $C$ and $E$.}\label{fig:Lface}
\end{center}
\end{figure}

Consider an example sensor system with 6 variables (columns) and 7 sensors
(rows) as summarized in the Table \ref{tab:7x6}. We will now show an example
categorification for the variable $L$ for violence using the pipeline
described in Section \ref{sec:bigpicture} and the machinery built up in this
report. This is a boolean variable with native data category \BOOL, and the
sensors that contribute to it are $C$, transit cams, and $E$, the Seattle
Times newspaper. We model this as an abstract simplicial complex with two
vertices ($C$ and $E$) and an edge ($L$), as shown in Figure \ref{fig:Lface}.
\begin{table}[h]
\begin{center}
  \begin{tabular}{l||c|c|c|c|c|c}
                    & $S$        & $O$           & $P$          & $I$        & $L$        & $R$        \\
                    & crowd Size & tOpic         & Place        & Intensity  & vioLence   & Role       \\
                    & Number     & Ontology term & Intersection & Level      & T/F        & Name       \\
                    & Scalar     & Partial order & Categorical  & Ordinal    & Boolean    & Categorical\\ \hline\hline
$A=$ police scAnner & \checkmark &               & \checkmark   &            &            & \checkmark \\ \hline
$C=$ transit Cams   &            &               & \checkmark   &            & \checkmark &            \\ \hline
$E=$ sEattle times  &            &               &              &            & \checkmark & \checkmark \\ \hline
$K=$ Komo news      & \checkmark &               & \checkmark   & \checkmark &            &            \\ \hline
$T_1=$ Twitter 1    &            & \checkmark    &              &            &            &            \\ \hline
$T_2=$ Twitter 2    & \checkmark & \checkmark    &              & \checkmark &            &            \\ \hline
$V=$ overhead Video & \checkmark &               &              &            &            &
  \end{tabular}
  \caption{The $7\times 6$ example.}\label{tab:7x6}
\end{center}
\end{table}

\begin{figure}[h]
\begin{center}
  \includegraphics{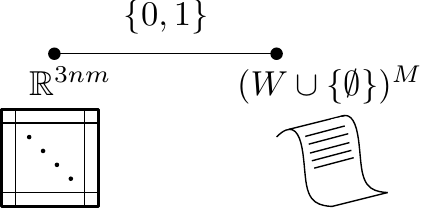}
  \caption{Raw data types for our sheaf.}\label{fig:Lface-raw}
\end{center}
\end{figure}

The zeroth step in this example, before we can define our analytics, is to
determine what our raw data feeds are. For the transit cameras let us assume
that they take static images (as opposed to video) which are $n\times m$
pixels. Then the raw data space for sensor $C$ would be $\R^{3n m}$ where
entries correspond to 3 color channels (red, green, blue) for each of the
$n\cdot m$ pixels.
The raw data space for the Seattle Times, sensor $E$, will be articles. Let
us assume that the Seattle Times has a word limit for each article, $M$, and
all articles are in English with word set $W$. Then a single article would be
an element of $\left(W\cup \{\emptyset\}\right)^M$, it is a vector of words
of length $M$ where the empty word is allowed (in case the article isn't
exactly length $M$). 
The simplicial complex with raw data types identified is shown in Figure
\ref{fig:Lface-raw}. The $\{0,1\}$ over the edge indicates the data space for
variable $L$.

These are simply the \emph{raw} data spaces for each data type. Now to
perform step 1, mathematization, we define the analytics on each sensor feed
for the variable $L$. For sensor $C$ we have two analytics, one for $P$
(Place) and one for $L$ (vioLence). Let us consider the analytic for $L$ to
be an image classification pipeline to determine probability of violence. We
can consider the target space of the analytic as the following object in \RV:
\begin{align*}
(\Omega, \mathcal{F}, P) &= (\R^{3nm}, \textit{Borel sets}, \textit{Image probabilities})\\
(S, \Sigma) &= ([0,1], \textit{Borel sets})\\
X &: \R^{3nm}\rightarrow [0,1]
\end{align*}
where the function $X$ is the result of an image classification algorithm
that takes in images and returns a probability, or confidence, that the image
contains violence. The ``Image probabilities'' would be a probability
distribution over all of the images, but this should not come into play in
the sheaf. We may just assume it is some probability distribution over the
set of all possible images. Recall in our discussion of the \RV{} category an
assignment is simply a choice of $\omega \in \Omega$ and an $s \in S$ such
that $X(\omega)=s$. Therefore, we define the analytic $f_{C,L}(\omega) =
X(\omega)$. 

Next we must define an analytic on $E$ (sEattle times) to inform $L$. Note
that there is an additional analytic for $R$ (Role) which we will not
consider. The analytic to inform $L$ from Seattle Times articles will be a
bag of words model. We first map $\left(W\cup \{\emptyset\}\right)^M$ to
$\N^{|W|+1}$ where each word vector is mapped to its vector of word (and
empty word) occurrence counts. Then, we can further select a set $V \subseteq
W$ of \emph{violent words} and a disjoint set $N\subseteq W$ of
\emph{non-violent} (or \emph{calm}) \emph{words} and project the space
$\N^{|W|+1}$ into $\N^{|V|+|N|}$ in the obvious fashion. Finally, we can map
further into $\N^2$ by summing up all violent word occurrences and separately
all non-violent word occurrences. So the full analytic is defined as
$f_{E,L}(a) = [\#\text{ of violent words in article }a, \#\text{ of calm
words in article }a]$.

Given our two analytics and target spaces in \RV{} and $\N^2$ we can do step
2 of our categorification pipeline and map both target spaces into the
category \BOOL, the native category for variable $L$. 
We define
\begin{align*}
f_1 &: \{(\omega,s) : X(\omega)=s\} \rightarrow \{0,1\}\\
f_2 &: \N^2 \rightarrow \{0,1\}
\end{align*}
where
\begin{align*}
f_1((\omega, s)) = \left\{\begin{array}{ll}
                          0 & s < 0.5 \\
                          1 & s \geq 0.5
                          \end{array}\right., \qquad
f_2((i,j) = \left\{\begin{array}{ll}
                    0 & i < j \\
                    1 & i \geq j
                    \end{array}\right..
\end{align*}
These two maps take the mathematized data space to the object $\{0,1\} \in
Ob(\BOOL)$ so that each element of the domain maps to an element of
$\{0,1\}$. Clearly these are not one-to-one functions, but that is not
required. Finally, we perform step 3 using the mapping defined in the
previous section to take $\{0,1\}$ to the vector space $\R[\{0,1\}]$.
%


\bibliographystyle{plain}
\bibliography{categorification}

\begin{thebibliography}{1}

\bibitem{BaMWeC2005}
Michael Barr and Charles Wells.
\newblock {\em Toposes, Triples and Theories}.
\newblock Number~12 in Reprints in Theory and Applications of Categories. 2005.

\bibitem{GoR2006}
Robert Goldblatt.
\newblock {\em {Topoi: The Categorical Analysis of Logic}}.
\newblock Dover, 2006.

\bibitem{RoM2016}
Michael Robinson.
\newblock Sheaves are the canonical datastructure for sensor integration.
\newblock arXiv:1603.01446v1 [math.AT].

\end{thebibliography}

\end{document}